\documentclass[notitlepage,12pt]{article}
\usepackage{amssymb,mathrsfs,amsmath,enumitem,comment}
\catcode`\@=11
\@addtoreset{equation}{section}

\catcode`\@=12

\newtheorem{Theorem}{Theorem}[section]
\newtheorem{Lemma}[Theorem]{Lemma}

\newtheorem{Remark}[Theorem]{Remark}

\def\R{{\mathbb R}}
\def\S{{\mathbb S}}
\def\C{C}
\def\Proof{\noindent\textit{Proof. }}
\def\qed{$~\square$\goodbreak \medskip}
\def\abstracts#1#2#3{{
	\centering{\begin{minipage}{5.0in}
	\small\baselineskip=13pt
	\parindent=0pt #1\\ \\
	\parindent=0pt #2\\ \\
	\parindent=0pt #3
\end{minipage}}\par}}

\begin{document}

\title{Radial graphs over domains of $\S^{n}$\\ with prescribed mean curvature}
\author{Paolo Caldiroli\footnote{The first author is partially supported by the PRIN2009 grant "Critical Point Theory and Perturbative Methods for Nonlinear Differential Equations".\hfill\eject Email of the authors: {\tt paolo.caldiroli@unito.it}, {\tt giovanni-gullo@hotmail.it}}~ and Giovanni Gullino\\ \normalsize{Dipartimento di Matematica, Universit\`a di Torino}\\
\normalsize{via Carlo Alberto, 10 -- 10123 Torino, Italy}}
\date{}
\maketitle

\abstracts{\textbf{Abstract.}
We prove the existence and uniqueness of radial graphs over a given domain of $\S^{n}$ having boundary on the sphere $\S^{n}$ and whose mean curvature at every point equals a prescribed positive function satisfying suitable barrier-type and monotonicity conditions.}

\section{Introduction and main result}
\label{S:Intro}

In this paper we are interested in hypersurfaces in $\R^{n+1}$, characterized as radial graphs over a given domain of $\S^{n}$, having boundary contained in $\S^{n}$ and whose mean curvature is a prescribed function $H\colon\R^{n+1}\to\R$. 

By \emph{radial graph over a domain $\Omega\subset\S^{n}$} we mean a hypersurface in $\R^{n+1}$ of the form
\begin{equation}
\label{eq:graph}
\Sigma=\Sigma(u)=\{e^{u(q)}q~|~q\in\overline\Omega\}
\end{equation}
for some mapping $u\in C^{2}(\Omega)\cap C^{0}(\overline\Omega)$. When $\Omega\ne\S^{n}$, we call \emph{$H$-bump on the sphere $\S^{n}$ supported by $\Omega$} a radial graph over $\Omega$ whose mean curvature at every point equals $H$ and such that $\partial\Sigma=\partial\Omega$. 

The problem of radial graphs over a given domain $\Omega$ in $\S^{n}$ has been studied firstly by J.~Serrin \cite{Se69} who proved existence and uniqueness for \emph{non-positive} prescribed mean curvature functions satisfying a suitable smallness condition involving the boundary datum and the geodesic mean curvature of $\partial\Omega$ with respect to $\S^{n}$. A key assumption is that the domain $\Omega$ is contained in an open hemisphere. 

Some variants, improvements, and extensions of the Serrin result in different directions have been later presented, even in recent years (see, e.g., \cite{AlDa07}, \cite{FuRi03}, \cite{Lo03}, \cite{LoMo99} and the references therein) but to our knowledge most of them concern radial graphs with non-positive \emph{constant} mean curvature. 

Here we focus on the case of \emph{prescribed} mean curvature and we discuss an existence result which in some sense is complementary to the Serrin one. In fact we deal with a situation in which $\Omega$ is \emph{any} strict regular domain in $\S^{n}$ with no restriction on its size neither on the curvature of its boundary, and $H$ is a \emph{positive} prescribed mean curvature function $H$ satisfying suitable barrier-type and monotonicity conditions. More precisely, we prove the following result.

\begin{Theorem}
\label{T:main}
Let $\Omega$ be a proper domain of $\S^{n}$ with boundary of class $C^{2,\alpha_{0}}$, for some $\alpha_{0}\in(0,1)$, and let
$$
A=\{\rho q\in\R^{n+1}~|~q\in\overline\Omega,~r_{1}\le\rho\le r_{2}\}
$$
with $0<r_{1}\le 1\le r_{2}<\infty$. Assume that $H$ is a real-valued function of class $C^{1}(A)$ satisfying:
\begin{gather}
\label{eq:barriers}
H(r_{1}q)\ge r_{1}^{-1}\quad\textrm{and}\quad H(r_{2}q)\le r_{2}^{-1}\quad\textrm{for every $q\in\overline\Omega$}\\
\label{eq:monotonicity}
\frac{\partial}{\partial\rho}\rho H(\rho q)\le 0\textrm{ for every $q\in\Omega$ and $\rho\in(r_{1},r_{2})$.}
\end{gather}
Then there exists an $H$-bump $\Sigma$ on the sphere $\S^{n}$ supported by $\Omega$. Moreover $\Sigma\subset A$ and it is the only $H$-bump supported by $\Omega$ and contained in $A$.
\end{Theorem}

We point out that our result is meaningful just for variable $H$. Indeed in case of constant mean curvature, the conditions (\ref{eq:barriers}) and (\ref{eq:monotonicity}) force the choice $r_{1}=r_{2}=1$ and $H\equiv 1$. In this case the result is trivial: the domain $\Omega$ itself is the hypersurface we look for. 

Notice that when $H\equiv 1$ and $\Omega$ is a spherical cap, there are two $H$-bumps supported by $\Omega$: the spherical cap $\Omega$ itself and that one obtained by reflecting $\S^{n}\setminus\Omega$ with respect to the hyperplane containing $\partial\Omega$. This lack of uniqueness does not violate the last statement in Theorem \ref{T:main} which in fact says that uniqueness holds in $A$, and in the previous case $A=\overline\Omega$. 

Actually the class of admissible mean curvature functions considered in our investigation is shaped on the mapping $H(X)=|X|^{-1}$. For such a function, conditions (\ref{eq:barriers}) and (\ref{eq:monotonicity}) hold true for any choice of $\Omega$, $r_{1}$ and $r_{2}$. We point out that the case $H(X)=|X|^{-1}$ exhibits dilation invariance but not invariance with respect to translation, as it happens in the constant case. Hence, Theorem \ref{T:main} can be viewed as a result about not necessarily too small perturbations of the mean curvature $|X|^{-1}$. In this direction see also \cite{Mu05} for a somehow similar problem. 

We also remark that the conditions required on the mean curvature function in Theorem \ref{T:main} were introduced in the papers \cite{BaKa74} and \cite{TrWe83} concerning the problem of hyperspheres with prescribed mean curvature (on this problem see also \cite{Ge98} and \cite{Tr85} and, for a different approach, concerning the parametric case \cite{CaMu04}, \cite{CaMu11} and \cite{Mu05}). 

Let us spend few words about the proof, which is accomplished by using classical tools of the Schauder theory like: a priori global estimates, the method of continuity, the Leray Schauder fixed point theorem. In fact we follow \emph{almost} the same argument developed in \cite{TrWe83}. Actually the proof sketched in \cite{TrWe83} contains a  misleading step (see Remark \ref{R:mistake}) regarding a technical feature, but it is possible to rectify it in a suitable way (see \cite{Gu12}). For this reason, even if the argument should be rather standard and known to experts, we preferred to display the proof in detail, both for the reader convenience and for better explain how to redress the misleading procedure suggested in \cite{TrWe83}.
\smallskip

\noindent
\textbf{Caution!} In many steps we refer to classical results concerning Dirichlet problems of the form 
\begin{equation}
\label{eq:flat}
\left\{\begin{array}{ll}\mathscr{L}u=f&\textrm{in $U$}\\ u=0&\textrm{on $\partial U$}\end{array}\right.
\end{equation}
where $U$ is a bounded domain in the flat Euclidean space $\R^{n}$ and $\mathscr{L}$ is a uniformly elliptic operator in $U$. We apply these results for Dirichlet problems of the form 
\begin{equation}
\label{eq:spherical}
\left\{\begin{array}{ll}Lu=f&\textrm{in $\Omega$}\\ u=0&\textrm{on $\partial\Omega$}\end{array}\right.
\end{equation}
where, in our context, $\Omega$ is a bounded domain in $\S^{n}$ and $L$ is a uniformly elliptic operator in $\Omega$. In doing this, by tacit agreement, we convert problem (\ref{eq:spherical}) into (\ref{eq:flat}) by stereographic projection of the sphere $\S^{n}$ into $\R^{n}$ from a point $P\in\S^{n}\setminus\overline\Omega$, playing the role of North pole, onto the Equatorial hyperplane, isomorphic to $\R^{n}$. In this way $\Omega$ is projected onto a bounded domain $U$ in $\R^{n}$ and the operator $L$ is converted by stereographic coordinates into a uniformly elliptic operator $\mathscr{L}$ in $U$. Hence the application to problem (\ref{eq:spherical}) of theorems holding for problem (\ref{eq:flat}) is justified.

\section{Proof}

\subsection{Proof of the existence with a stronger assumption}
\label{SS:1}
Firstly we prove the theorem under a stronger condition on $H$ and precisely:
\begin{equation}
\label{eq:strict-barriers}
H(r_{1}q)>r_{1}^{-1}\quad\textrm{and}\quad H(r_{2}q)<r_{2}^{-1}\quad\textrm{for every $q\in\overline\Omega$}\end{equation}
instead of (\ref{eq:barriers}). Note that (\ref{eq:strict-barriers}) implies that $r_{1}<r_{2}$. We extend $H$ to a mapping still denoted $H$ on the cone $\widehat{\Omega}=\{\rho q~|~q\in\overline\Omega,~\rho>0\}$ in the following way: for every $q\in\overline\Omega$ set
$$
h_{1}(q)=\left[\frac{\partial}{\partial\rho}\rho H(\rho q)\right]_{\rho=r_{1}}\quad\textrm{and}\quad h_{2}(q)=\left[\frac{\partial}{\partial\rho}\rho H(\rho q)\right]_{\rho=r_{2}}
$$
and
\begin{equation}
\label{eq:extension}
{H}(\rho q)=\left\{\begin{array}{ll}
\dfrac{r_{1}}{\rho}~\!H(r_{1}q)+\left(1-\dfrac{r_{1}}{\rho}\right)h_{1}(q)&\textrm{for $\rho\in(0,r_{1})$}\vspace{4pt}\\
H(\rho q)&\textrm{for $\rho\in[r_{1},r_{2}]$}\vspace{4pt}\\
\dfrac{r_{2}}{\rho}~\!H(r_{2}q)+\left(1-\dfrac{r_{2}}{\rho}\right)h_{2}(q)&\textrm{for $\rho\in(r_{2},\infty)$}
\end{array}\right.
\end{equation}
One plainly checks that the mapping $H$ is of class $\C^{1}$ on $\widehat{\Omega}$ and satisfies
\begin{equation}
\label{eq:strict-barriers2}
\left\{\begin{array}{l}
H(X)>|X|^{-1}\textrm{ if $|X|\le r_{1},~X\in\widehat{\Omega}$}\\
H(X)<|X|^{-1}\textrm{ if $|X|\ge r_{2},~X\in\widehat{\Omega}$}.
\end{array}\right.
\end{equation}

The next step consists in stating the problem of radial $H$-graphs over domains in $\S^{n}$ satisfying some boundary condition in terms of a Dirichlet problem for a quasilinear elliptic equation. This is the content of what follows.

\begin{Lemma}
\label{L:dirichlet-problem}
Let $\Omega$ be a smooth domain in $\S^{n}$, $\Omega\ne\S^{n}$. A radial graph on $\Omega$ is an $H$-bump if and only if the mapping $u\in C^{2}(\Omega)\cap C^{0}(\overline\Omega)$ defining the radial graph via (\ref{eq:graph}) satisfies the Dirichlet problem
\begin{equation}
\label{eq:dirichlet}
\left\{\begin{array}{ll}
\left((1+|\nabla u|^{2})\delta_{ij}-u_{i}u_{j}\right)u_{ij}& \\
\phantom{1+|\nabla u|^{2})}=n(1+|\nabla u|^{2})\left(1-\sqrt{1+|\nabla u|^{2}}~\!e^{u}~\!H(e^{u}q)\right)&\textrm{in }\Omega\\
u=0&\textrm{on }\partial\Omega\end{array}
\right.
\end{equation}
where the subscripts denote covariant derivatives in an orthogonal frame on $\S^{n}$, and $\nabla$ is the gradient operator with the standard metric of $\S^{n}$. 
\end{Lemma}

\Proof 
See \cite{TrWe83} or \cite{Lo03}.
\qed

Our next goal is to convert the Dirichlet problem (\ref{eq:dirichlet}) into a suitable fixed point equation in some functional space. For this reason we begin by studying the linear operator 
\begin{equation}
\label{eq:Lw}
L_{w}u=\left((1+|\nabla w|^{2})\delta_{ij}-w_{i}w_{j}\right)u_{ij}.
\end{equation}
where $w$ is a fixed, sufficiently regular mapping on $\Omega$. We work with the standard spaces $\C^{k}(\overline\Omega)$ and $\C^{k,\alpha}(\overline\Omega)$ endowed with their natural norms, denoted $\|\cdot\|_{k}$ and $\|\cdot\|_{k,\alpha}$, respectively ($k$ is a non-negative integer and $\alpha\in(0,1)$). Because of the boundary condition in problem (\ref{eq:dirichlet}) it is also convenient to introduce the spaces 
$$
\C^{k,\alpha}_{0}(\overline\Omega):=\{u\in \C^{k,\alpha}(\overline\Omega)~|~u|_{\partial\Omega}=0\}.
$$
Fixing $w\in\C^{1,\alpha}(\overline\Omega)$, the operator $L_{w}$ defined in (\ref{eq:Lw}) is a linear bounded operator from $\C^{2,\alpha}_{0}(\overline\Omega)$ in $\C^{0,\alpha}(\overline\Omega)$ because $w\in\C^{1}(\overline\Omega)$. In fact the following result holds.

\begin{Lemma}
\label{L:Lw}
If $\alpha\in(0,\alpha_{0}]$ then for every $w\in\C^{1,\alpha}(\overline\Omega)$ the operator $L_{w}$ is a bijection of $\C^{2,\alpha}_{0}(\overline\Omega)$ onto $\C^{0,\alpha}(\overline\Omega)$.
\end{Lemma}

\Proof
A mapping $u\in\C^{2,\alpha}_{0}(\overline\Omega)$ belongs to the kernel of $L_{w}$ if and only if $u$ solves the Dirichlet problem
$$
\left\{\begin{array}{ll}
a^{ij}u_{ij}=0&\textrm{in }\Omega\\
u=0&\textrm{on }\partial\Omega\end{array}
\right.\quad\textrm{where}\quad a^{ij}=(1+|\nabla w|^{2})\delta_{ij}-w_{i}w_{j}.
$$
Since $a^{ij}\xi_{i}\xi_{j}=(1+|\nabla w|^{2})|\xi|^{2}-\langle\nabla w,\xi\rangle^{2}$ we readily obtain that 
\begin{equation}
\label{eq:ellipticity}
|\xi|^{2}\le a^{ij}\xi_{i}\xi_{j}\le\left(1+2\sup_{\Omega}|\nabla w|^{2}\right)|\xi|^{2}
\end{equation}
that is, the operator $L_{w}u=a^{ij}u_{ij}$ is uniformly elliptic. 
Then, by the maximum principle, since $u=0$ on $\partial\Omega$, we infer that $u=0$ in $\Omega$. This shows that $L_{w}$ is injective.   In order to prove that $L_{w}$ is onto, we use the continuity method (see Section 5.2 in \cite{GiTr98}). More precisely we introduce a family of operators $\mathscr{L}_{t}\colon\C^{2,\alpha}_{0}(\overline\Omega)\to\C^{0,\alpha}(\overline\Omega)$ with $t\in[0,1]$ defined by
$$
\mathscr{L}_{t}=(1-t)L_{0}+tL_{w}~\!.
$$
Notice that $\mathscr{L}_{0}=L_{0}$ is the Laplace-Beltrami operator on $\S^{n}$ and for every $f\in\C^{0,\alpha}(\overline\Omega)$ the Dirichlet problem 
$$
\left\{\begin{array}{ll}
\Delta_{\S^{n}}u=f&\textrm{in }\Omega\\
u=0&\textrm{on }\partial\Omega\end{array}
\right.
$$
admits a solution $u\in\C^{2,\alpha}(\overline\Omega)$. That is, $\mathscr{L}_{0}$ sends $\C^{0,\alpha}(\overline\Omega)$ onto $\C^{2,\alpha}_{0}(\overline\Omega)$. Now we claim that there exists a constant $C>0$ such that
\begin{equation}
\label{eq:continuity}
\|u\|_{2,\alpha}\le C\|\mathscr{L}_{t}u\|_{0,\alpha}\quad\textrm{for every $t\in[0,1]$ and for every $u\in\C^{2,\alpha}_{0}(\overline\Omega)$.}
\end{equation}
In view of the method of continuity (Theorem 5.2 in \cite{GiTr98}) this is enough to infer that $L_{w}=\mathscr{L}_{1}$ is onto.
We show (\ref{eq:continuity}) arguing by contradiction. If (\ref{eq:continuity}) is false then there exist sequences $\{t_{k}\}\subset[0,1]$, $\{u_{k}\}\subset\C^{2,\alpha}_{0}(\overline\Omega)$ such that
\begin{equation}
\label{eq:contr1}
\|\mathscr{L}_{t_{k}}u_{k}\|_{0,\alpha}\to 0\quad\textrm{and}\quad\|u_{k}\|_{2,\alpha}=1.
\end{equation}
By compactness, in particular using also the Ascoli-Arzel\`a Theorem, there exist $t\in[0,1]$ and $u\in\C^{2}_{0}(\overline\Omega)$ such that, up to subsequences,
$$
t_{k}\to t\quad\textrm{and}\quad u_{k}\to u\textrm{ in }\C^{2}(\overline{\Omega}).
$$
By continuity we get $\mathscr{L}_{t}u=0$. Since $\mathscr{L}_{t}$ is a convex combination of elliptic operators, it is so, too. Hence $u=0$. In particular
\begin{equation}
\label{eq:contr2}
u_{k}\to 0\textrm{ in }\C^{0}(\overline{\Omega}).
\end{equation}
We point out that 
$$
\mathscr{L}_{t}u=a^{ij}_{t}u_{ij}\quad\textrm{where}\quad a^{ij}_{t}=\left((1+t|\nabla w|^{2})\delta_{ij}-tw_{i}w_{j}\right)u_{ij}
$$
and, as in (\ref{eq:ellipticity}), 
$$
|\xi|^{2}\le a^{ij}_{t}\xi_{i}\xi_{j}\le\left(1+2\sup_{\Omega}|\nabla w|^{2}\right)|\xi|^{2}\quad\textrm{for every }t\in[0,1].
$$
Notice that the ellipticity constants are independent of $t$. Taking into account of the $C^{2,\alpha_{0}}$ regularity of the domain, 
we can apply the a priori global Schauder estimates (Theorem 6.6 in \cite{GiTr98}) 
obtaining that
$$
\|u_{k}\|_{2,\alpha}\le C\left(\|u_{k}\|_{0}+\|\mathscr{L}_{t_{k}}u_{k}\|_{0,\alpha}\right)
$$
with $C$ independent of $k$. This yields a contradiction with (\ref{eq:contr1}) and (\ref{eq:contr2}). Hence (\ref{eq:continuity}) holds true and the proof is complete.
\qed

\noindent
From now on we take $\alpha\in(0,\alpha_{0}]$.

\begin{Lemma} [uniform bound on the operator norms] 
\label{L:bounds}
For every $C>0$ there exists $K>0$ such that if $\|w\|_{1,\alpha}\le C$ then $\|u\|_{2,\alpha}\le K\|L_{w}u\|_{0,\alpha}$ for every $u\in\C^{2,\alpha}_{0}(\overline\Omega)$.
\end{Lemma}

\Proof
We argue by contradiction as in the last part of the proof of Lemma \ref{L:Lw}. If the result is false then there exist a bounded sequence $\{w_{k}\}$ in $\C^{1,\alpha}(\overline\Omega)$ and a sequence $\{u_{k}\}$ in $\C^{2,\alpha}_{0}(\overline\Omega)$ such that 
\begin{equation}
\label{eq:contr3}
\|u_{k}\|_{2,\alpha}=1\quad\textrm{and}\quad\|L_{w_{k}}u_{k}\|_{0,\alpha}\to 0.
\end{equation}
By compactness, there exist $w\in\C^{1}(\overline\Omega)$ and $u\in\C^{2}_{0}(\overline\Omega)$ such that, up to subsequences,
$$
w_{k}\to w\textrm{ in }\C^{1}(\overline{\Omega})\quad\textrm{and}\quad u_{k}\to u\textrm{ in }\C^{2}(\overline{\Omega}).
$$
By continuity we get $L_{w}u=0$. Then $u=0$, by Lemma \ref{L:Lw}. Taking into account of (\ref{eq:ellipticity}), we observe that the operators $L_{w_{k}}$ are uniformly elliptic with ellipticity constants independent of $k$ (but depending on $C$). Using the a priori global Schauder estimates (Theorem 6.6 in \cite{GiTr98}),  
we obtain that
$$
\|u_{k}\|_{2,\alpha}\le\widetilde{C}\left(\|u_{k}\|_{0}+\|L_{w_{k}}u_{k}\|_{0,\alpha}\right)
$$
with $\widetilde{C}$ independent of $k$. Since $u_{k}\to 0$ in $\C^{0}(\overline{\Omega})$ and by (\ref{eq:contr3}) we reach a contradiction.~\qed

\begin{Lemma} 
\label{L:compactness}
Let $\{w_{k}\}$ be a bounded sequence in $\C^{1,\alpha}(\overline\Omega)$ and let $\{f_{k}\}$ be a bounded sequence in $\C^{0,\alpha}(\overline\Omega)$. Then the sequence $\{u_{k}\}$ of solutions of 
\begin{equation}
\label{eq:k}
\left\{\begin{array}{ll}L_{w_{k}}u_{k}=f_{k}&\textrm{in }\Omega\\ u_{k}=0&\textrm{on }\partial\Omega\end{array}\right.
\end{equation}
is bounded in $\C^{2,\alpha}(\overline\Omega)$.
\end{Lemma}

\Proof
The existence of $u_{k}$ is guaranteed by Lemma \ref{L:Lw}. The conclusion plainly follows from Lemma \ref{L:bounds}.
\qed

Let us introduce a family of operators $T_{t}\colon\C^{1,\alpha}(\overline\Omega)\to\C^{1,\alpha}(\overline\Omega)$ depending on a parameter $t\in[0,1]$, defined as follows: for every $w\in\C^{1,\alpha}(\overline\Omega)$ let $T_{t}w=u$ be the unique solution of problem
$$
\left\{\begin{array}{ll}L_{w}u=nt(1+|\nabla w|^{2})\left(1-\sqrt{1+|\nabla w|^{2}}~\!e^{w}~\!H(e^{w}q)\right)&\textrm{in }\Omega\\ u=0&\textrm{on }\partial\Omega.\end{array}\right.
$$
We point out that the operator $T_{t}$ is well defined in view of Lemma \ref{L:Lw} and in fact takes values in $\C^{2,\alpha}_{0}(\overline\Omega)$. Moreover one plainly recognizes that $T_{t}=tT_{1}$ and that fixed points of $T_{1}$ are solutions of problem (\ref{eq:dirichlet}). 

\begin{Remark}
\label{R:mistake}
The definition of the operators $T_{t}$ considered here differs from that one in \cite{TrWe83}. In that paper the authors define $T_{t}w$ for $w\in\C^{1,\alpha}(\S^{n})$ as the unique solution $u\in\C^{2,\alpha}(\S^{n})$ of
$$
\mathrm{div}_{\S^{n}}\left((1+|\nabla w|^{2})^{-\frac{1}{2}}\nabla u\right)-u=t\left(n(1+|\nabla w|^{2})^{-\frac{1}{2}}-n~\!e^{w}~\!H(e^{w}q)-w\right).
$$
But the differential operator at the left hand side is an operator in divergence form with coefficients which are just in $\C^{0,\alpha}$ and it is well known that this is not enough to guarantee the existence of $\C^{2,\alpha}$ solutions. A way to redress the argument and recover the $C^{0}$ and $C^{1}$ estimates proved in \cite{TrWe83} is to introduce the operators $L_{w}$ given by (\ref{eq:Lw}). These operators are invertible and the reasoning can be carried through without trouble. In a similar way one can treat also the problem studied in \cite{TrWe83} (see \cite{Gu12}).
\end{Remark}

In order to find that $T_{1}$ does possess a fixed point we will apply the Leray-Schauder Theorem (Theorem 11.3 in \cite{GiTr98}). To this extent we begin checking the following compactness property.

\begin{Lemma} 
\label{L:T-compact}
The operator $T_{1}$ is compact in $\C^{1,\alpha}(\overline\Omega)$ for every $\alpha\in(0,1)$.
\end{Lemma}

\Proof
Let $\{w_{k}\}$ be a bounded sequence in $\C^{1,\alpha}(\overline\Omega)$ and set
$$
f_{k}=n(1+|\nabla w_{k}|^{2})\left(1-\sqrt{1+|\nabla w_{k}|^{2}}~\!e^{w_{k}}~\!H(e^{w_{k}}q)\right).
$$
Then the sequence $\{f_{k}\}$ is bounded in $\C^{0,\alpha}(\overline\Omega)$. The mapping $u_{k}=T_{1}w_{k}$ is the unique solution of problem (\ref{eq:k}). Then the conclusion follows by means of Lemma \ref{L:compactness} and by the Ascoli-Arzel\`a Theorem.
\qed

We also need a priori $C^{1,\alpha}$ uniform estimates on fixed points of $T_{t}$, for all $t\in[0,1]$. Following a standard scheme, we start by finding $C^{0}$ estimates. To this goal we use the barrier condition (\ref{eq:barriers}).

\begin{Lemma} [A priori $C^{0}$ estimates]
\label{L:C0-bound}
For every $t\in[0,1]$ if $u$ is a fixed point of $T_{t}$ then 
$\log r_{1}\le u(q)\le\log r_{2}$ for every $q\in\overline\Omega$.
\end{Lemma}

\Proof
Let $t\in[0,1]$ and let $u$ be a fixed point of $T_{t}$. Then $u$ is a $C^{2}$ solution of the Dirichlet problem
$$
\left\{\begin{array}{ll}L_{u}u=nt(1+|\nabla u|^{2})\left(1-\sqrt{1+|\nabla u|^{2}}~\!e^{u}~\!H(e^{u}q)\right)&\textrm{in }\Omega\\ u=0&\textrm{on }\partial\Omega\end{array}\right.
$$
Let $q_{0}\in\overline\Omega$ be such that $u(q_{0})=\max_{\overline\Omega}u$. Assume by contradiction $u(q_{0})>\log r_{2}$. Then 
$q_{0}\in\Omega$ because $r_{2}\ge 1$ and $u=0$ on $\partial\Omega$, and then $\nabla u(q_{0})=0$ and
$$
L_{u}u(q_{0})\le 0.
$$
On the other hand it must be $t>0$ because otherwise $u=0$, by Lemma \ref{L:Lw}. Moreover
$$
L_{u}u(q_{0})=nte^{u(q_{0})}\left(\frac{1}{e^{u(q_{0})}}-H(e^{u(q_{0})}q_{0})\right)>0
$$
because $H(X)<|X|^{-1}$ as $|X|>r_{2}$, see (\ref{eq:strict-barriers2}). Thus we reach a contradiction. The same argument holds in order to show that $\min_{\Omega}u\ge \log r_{1}$.
\qed

Now we discuss the a priori $C^{1,\alpha}$ uniform bound on fixed points of $T_{t}$. As we will see, the monotonicity assumption (\ref{eq:monotonicity}) enters in order to get gradient estimates.

\begin{Lemma} [A priori $C^{1,\alpha}$ estimates]
\label{L:bound}
There exists $\alpha\in(0,\alpha_{0}]$ and $C>0$ such that for every $t\in[0,1]$, if $u$ is a fixed point of $T_{t}$ then $\|u\|_{1,\alpha}\le C$. 
\end{Lemma}

\Proof
We just give a sketch since the proof is essentially the same as in \cite{TrWe83}. One uses a gradient estimate proved by Treibergs and Wei \cite{TrWe83} for a class of equations shaped on the equation of prescribed mean curvature. More precisely it deals with equations of the form
$$
a^{ij}u_{ij}=b\left(q,u,|\nabla u|^{2}\right)
$$
where $a^{ij}=\left(1+|\nabla u|^{2}\right)\delta_{ij}-u_{i}u_{j}$, and $b=b(q,u,v)\in C^{1}(\Omega\times\R\times[0,\infty))$ satisfying
\begin{equation}
\label{eq:TW}
|b_{q}|\le C_{1}(1+v)^{\frac{3}{2}},\quad b_{u}\ge -C_{2}(1+v)~\!,\quad (1+3v)b-2(1+v)vb_{v}\ge -C_{3}(1+v)^{\frac{3}{2}}
\end{equation}
for some nonnegative constants $C_{1},C_{2},C_{3}$. We are in this setting with 
$$
b(q,u,v)=nt(1+v)\left(1-\sqrt{1+v}~\!e^{u}~\!H(e^{u}q)\right)
$$
and conditions (\ref{eq:TW}) are satisfied by taking
\begin{equation}
\label{eq:TW-constants}
C_{1}=\max_{X\in A}|X|^{2}|\nabla H(X)|~\!,\quad C_{2}=0~\!,\quad C_{3}=\max_{X\in A}|X|H(X)~\!.
\end{equation}
As in \cite{TrWe83} we can conclude that there exist constants $C_{4}$ and $C_{5}$ such that for every $t\in[0,1]$, if $u$ is a fixed point of $T_{t}$ then
\begin{equation}
\label{eq:C1-bound}
|\nabla u|^{2}\le C_{4}(r_{2}/r_{1})^{C_{5}}.
\end{equation}
Having found uniform bounds on the $C^{0}$ and $C^{1}$ norms of $u$, in order to get a priori estimates in $C^{1,\alpha}$ for some $\alpha\in(0,\alpha_{0}]$, one uses a standard result (Theorem 13.2 in  \cite{GiTr98}).~\qed

By Lemmas \ref{L:T-compact} and \ref{L:bound} the Leray-Schauder theorem can be applied to the fixed point equation $u=tT_{1}u$ in $\C^{1,\alpha}(\overline\Omega)$ and one obtains the existence of a fixed point for $T_{1}$, namely a solution $u$ of (\ref{eq:dirichlet}). Then, by Lemma \ref{L:dirichlet-problem}, 
the radial graph $\Sigma(u)$ defined by (\ref{eq:graph}) is an $H$-bump on the sphere $\S^{n}$ supported by $\Omega$. This completes the proof of existence when $H\in\C^{1}(A)$ satisfies (\ref{eq:monotonicity}) and (\ref{eq:strict-barriers}).
\qed

\subsection{Proof of the existence without the extra assumption (\ref{eq:strict-barriers})}

Assume that $r_{1}<r_{2}$ and that $H\in\C^{1}(A)$ verifies (\ref{eq:barriers}) and (\ref{eq:monotonicity}). Define $H_{\varepsilon}(X)=|X|^{-\varepsilon}H(X)$, with $\varepsilon>0$ and $X\in A$. Notice that $H_{\varepsilon}\in\C^{1}(A)$ satisfies (\ref{eq:monotonicity}) and (\ref{eq:strict-barriers}). Then by what proved in Subsection \ref{SS:1}, there exists $u_{\varepsilon}\in \C^{2,\alpha}(\overline\Omega)$ solving problem (\ref{eq:dirichlet}) with $H_{\varepsilon}$ instead of $H$. By Lemma \ref{L:C0-bound} we have that $\log r_{1}\le u_{\varepsilon}(q)\le\log r_{2}$ for every $q\in\overline\Omega$ and for every $\varepsilon>0$. Moreover the constants $C_{4}$ and $C_{5}$ in the a priori gradient estimate (\ref{eq:C1-bound}) can be taken independent of $\varepsilon\in(0,1]$ because they depend just on $r_{1}$, $r_{2}$ and on the constants $C_{1},C_{2},C_{3}$ in (\ref{eq:TW-constants}) which can be bounded uniformly with respect to $\varepsilon$ since $H_{\varepsilon}\to H$ in $C^{1}(A)$. Then also a uniform bound for $\|u_{\varepsilon}\|_{1,\alpha}$ holds, thanks to Theorem 13.2 in \cite{GiTr98}. Using the mean curvature equation, one infer that $\|u_{\varepsilon}\|_{2,\alpha}\le C$ for every $\varepsilon>0$ small enough and for some constant  $C$ independent of $\varepsilon$. By the Ascoli-Arzel\`a Theorem there exists $u\in\C^{2}(\overline\Omega)$ such that $u_{\varepsilon}\to u$ in $\C^{2}(\overline\Omega)$ for a sequence $\varepsilon\to 0$. Passing to the limit in the mean curvature equation satisfied by $u_{\varepsilon}$ one concludes that $u$ is a solution of problem (\ref{eq:dirichlet}). The case $r_{1}=r_{2}=1$ is trivial, as observed in Section \ref{S:Intro}. 
\qed

\subsection{Proof of the uniqueness}
Let $u,v\in\C^{2}(\Omega)\cap C^{0}(\overline\Omega)$ be two solutions of problem (\ref{eq:dirichlet}) such that the corresponding radial graphs $\Sigma(u)$ and $\Sigma(v)$ are contained in $A$. Let us consider the extension of $H$ on the cone $\widehat{\Omega}=\{\rho q~|~q\in\overline\Omega,~\rho>0\}$ defined by (\ref{eq:extension}). Notice that 
\begin{equation}
\label{eq:monotonicity2}
\frac{\partial}{\partial\rho}\rho H(\rho q)\le 0\quad\textrm{for every $\rho>0$ and $q\in\Omega$.}
\end{equation}
If $u\ne v$ then there exists $\overline{q}\in\Omega$ such that $u(\overline{q})\ne v(\overline{q})$. We can assume $u(\overline{q})<v(\overline{q})$. Then there exists $\mu>0$ such that $u(q)+\mu\ge v(q)$ for every $q\in\Omega$ and $u(q_{0})+\mu=v(q_{0})$ at some $q_{0}\in\Omega$. Set $\widetilde{u}=u+\mu$ and observe that $\widetilde{u}$ solves
$$
\left((1+|\nabla\widetilde{u}|^{2})\delta_{ij}-\widetilde{u}_{i}\widetilde{u}_{j}\right)\widetilde{u}_{ij}\le n(1+|\nabla \widetilde{u}|^{2})\left(1-\sqrt{1+|\nabla \widetilde{u}|^{2}}~\!e^{\widetilde{u}}~\!H(e^{\widetilde{u}}q)\right)\textrm{ in $\Omega$}
$$
because of (\ref{eq:monotonicity2}) and since $\mu>0$. Notice that the radial graph defined by $\widetilde{u}$ stays over (in the radial direction) that one corresponding to $v$ and they touch themselves at the point $X_{0}=q_{0}e^{v(q_{0})}$. Now we compare $\widetilde{u}$ and $v$ by means of the Hopf maximum principle. In particular we use the following version, stated in \cite{PuSe04}, Theorem 2.3.

\begin{Theorem} [Touching Lemma]
\label{T:touching}
Let $u_{1}$ and $u_{2}$ be $C^{2}$ solutions of the nonlinear differential inequalities
$$
\mathscr{F}(x,u_{1},Du_{1},D^{2}u_{1})\ge 0~,\quad\mathscr{F}(x,u_{2},Du_{2},D^{2}u_{2})\le 0
$$
in a domain $U$ in $\R^{n}$, with $\mathscr{F}$ of class $C^{1}$ in the variables $u,Du,D^{2}u$. Suppose also that the matrix
$$
Q_{ij}=\frac{\partial\mathscr{F}}{\partial u_{ij}}(x,u_{1},Du_{1},\theta D^{2}u_{1}+(1-\theta)D^{2}u_{2})
$$
is positive definite in $U$ for all $\theta\in[0,1]$.
If $u_{1}\le u_{2}$ in $U$ and $u_{1}(x_{0})=u_{2}(x_{0})$ at some point $x_{0}\in U$, then $u_{1}=u_{2}$ in $U$.
\end{Theorem}
We are in position to apply Theorem \ref{T:touching} with
\begin{equation*}
\begin{split}
\mathscr{F}(x,u,Du,D^{2}u)=&\left((1+|\nabla u|^{2})\delta_{ij}-u_{i}u_{j}\right)u_{ij}\\
&\qquad-n(1+|\nabla u|^{2})\left(1-\sqrt{1+|\nabla u|^{2}}~\!e^{u}~\!H(e^{u}q)\right)
\end{split}
\end{equation*}
where $x$ is the stereographic projection of $q$ and the derivatives in the arguments of $\mathscr{F}$ are meant in stereographic coordinates (see the Caution note at the end of the Introduction). In our case $u_{1}=v$, $u_{2}=\widetilde{u}$ and $Q_{ij}=(1+|\nabla v|^{2})\delta_{ij}-v_{i}v_{j}$ which is a positive definite matrix thanks to (\ref{eq:ellipticity}). The assumptions of Theorem \ref{T:touching} are fulfilled and thus we deduce that $\widetilde{u}=v$ in $\Omega$. This is impossible since $\widetilde{u}|_{\partial\Omega}=\mu>0=v|_{\partial\Omega}$. Hence it must be $u=v$. 
\qed

\end{document}